\documentclass[fleqn,oneside]{article}
\usepackage{indentfirst}
\usepackage{verbatim}
\usepackage{fullpage}

\input{xy}
\xyoption{all}
\xyoption{matrix}
\CompileMatrices

\usepackage{amsmath,amsthm,amsfonts}
\theoremstyle{plain}
\newtheorem{teo}{Theorem}[section]

\newtheorem{lem}[teo]{Lemma}
\newtheorem{prop}[teo]{Proposition}

\theoremstyle{remark}
\newtheorem*{remark}{Remark}

\newcommand{\exs}{\noindent\textbf{Examples: }}

\def\CC{{\mathbb{C}}}
\def\ZZ{{\mathbb{Z}}}
\def\NN{{\mathbb{N}}}
\def\U{{\cal U}}
\def\Z{{\cal Z}}
\def\G{{\cal G}}
\def\Om{\Omega}

\def\Id{I\!d}
\def\HH{H\!H}
\def\PSL{\mathit{PSL}}
\def\SL{\mathit{SL}}
\def\Im{\mathrm{Im}}

\let\ot\otimes
\let\wt\widetilde
\def\cl#1{{\langle #1\rangle}}

\DeclareMathOperator{\Aut}{Aut}
\DeclareMathOperator{\gr}{gr}
\DeclareMathOperator{\ad}{ad}
\DeclareMathOperator{\Tor}{Tor}
\DeclareMathOperator{\Ker}{Ker}
\DeclareMathOperator{\Hom}{Hom}
\DeclareMathOperator{\chr}{char}

\title{Hochschild homology and cohomology of generalized Weyl 
algebras}

\author{Marco A. Farinati${}^1$ 
\and Andrea Solotar${}^{1,2}$
\and Mariano Su\'arez-\'Alvarez${}^1$}

\begin{document}
\maketitle
\tableofcontents
\footnotetext[1]
{Dto. de Matem\'atica, Facultad de Cs. Exactas y Naturales.
Universidad de Buenos Aires. Ciudad Universitaria Pab I. 1428,
Buenos Aires - Argentina. e-mail: \texttt{mfarinat@dm.uba.ar},
\texttt{asolotar@dm.uba.ar},  \texttt{mariano@dm.uba.ar}\\
Research partially supported by UBACYT TW69 and CONICET.}

\footnotetext[2]
{Research member of CONICET (Argentina).}

%%%%%%%%%%%%%%%%%%%%%%%%%%%%%%%%%%%%%%%%%%%%%%%%%%%%%%%%%%
\section*{Introduction}

The relevance of algebras such as the Weyl algebra $A_n(\CC)$, the
enveloping algebra $\U(\mathfrak{sl}_2)$ and its primitive quotients
$B_{\lambda}$ and other algebras related to algebras of differential operators
is already well-known.  Recently, some articles where their
Hochschild homology and cohomology has an important role have been written
(see for example \cite{afls},  \cite{al1}, \cite{al2},  \cite{gui-eti},
\cite{odile},  \cite{mariano}, \cite{marianococprim}).  Both the results
obtained in \cite{afls} and in \cite{marianococprim} seem to depend
strongly on intrinsic properties of $A_1(\CC)$ and $\U(\mathfrak{sl}_2)$.
However, it is not entirely the case.

In this article we consider a class of algebras, called generalized Weyl
algebras (GWA for short), defined by V.~Bavula in \cite{Bav} and studied by
himself and collaborators in a series of papers (see for example \cite{Bav},
\cite{BavJor}, \cite{BavLen}) from the point of view of ring theory.

Our aim is to compute the Hochschild homology and cohomology groups of these
algebras and to study whether there is a duality between these groups.

Examples of GWA are, as we said before, $n$-th Weyl algebras,
$\U(\mathfrak{sl}_2)$, primitive quotients of $\U(\mathfrak{sl}_2)$, and
also the subalgebras of invariants of these algebras under the action of
finite cyclic subgroups of automorphisms.

As a consequence we recover in a simple way the results of \cite{afls} and we
also complete results of \cite{odile},  \cite{michelinkassel} and of
\cite{marianococprim}, giving at the same time a unified method for GWA.

\medskip

The article is organized as follows:

In section \ref{sect:GWA} we recall from \cite{Bav} the definition of generalized Weyl
algebras and state the main theorems.

In section \ref{sect:res} we describe the resolution used afterwards in order
to compute the Hochschild homology and cohomology groups. We prove
a ``reduction" result (Proposition \ref{prop:reduction}) and finally we 
prove the main theorem for homology using an spectral sequence argument.

Section \ref{sect:coho} is devoted to the computation of Hochschild cohomology of
GWAs. As a consequence of the results we obtain, we notice that the hypotheses
of Theorem 1 of \cite{VdB} are not sufficient to assure duality between
Hochschild homology and cohomology. We state here hypotheses under
which duality holds.

In section \ref{sect:inv} we consider subalgebras of invariants 
of the previous ones under
diagonalizable cyclic actions of finite order. We first show that these
subalgebras are also GWA and we state the main theorem concerning subalgebras
of invariants.

Finally, in section \ref{sect:apps}, we begin by describing some applications of 
the above results. The first application is to specialize the results to
the usual Weyl algebra. Secondly, we consider the primitive quotients
of ${\frak sl}_2$, and considering the Cartan involution $\Om$ we answer
a question of Bavula (\cite{BavJor} remark 3.30) and we finish the
proof of the main theorem. The formula for the dimension of $\HH_*(A^G)$
explains, in particular, the computations made by O.~Fleury for
$\HH_0(B_{\lambda}^G)$.

We will work over a field $k$ of characteristic zero and all algebras will be
$k$-algebras. Given a $k$-algebra $A$, $\Aut_k(A)$ will always denote the
group of $k$-algebra automorphisms of $A$.

%%%%%%%%%%%%%%%%%%%%%%%%%%%%%%%%%%%%%%%%%%%%%%%%%%%%%%%%%%
\section{Generalized Weyl Algebras}\label{sect:GWA}

We recall the definition of generalized Weyl algebras given by Bavula in
\cite{Bav}.

Let $R$ be an algebra, fix a central element
$a\in\Z(R)$ and $\sigma\in \Aut_k(R)$. The
generalized Weyl Algebra $A=A(R,a,\sigma)$ is the $k$-algebra
generated by $R$ and two new free variables $x$ and $y$ subject
to the relations:
  \begin{alignat*}{2}
  yx&=a &\qquad&xy=\sigma(a) \\
\intertext{and}  
  xr&=\sigma(r)x &&ry=y\sigma(r)
  \end{alignat*}
for all $r\in R $.

\medskip

\exs\begin{enumerate}

\item If $R=k[h]$, $a=h$ and $\sigma\in \Aut_k(k[h])$ is the unique
automorphism determined by $\sigma(h)=h-1$; then $A(k[h],a,\sigma)\cong
A_1(k)$, the usual Weyl algebra, generated by $x$ and $y$ subject to the relation
$[x,y]=1$.

\item Let $R=k[h,c]$, $\sigma(h)=h-1$ and $\sigma(c)=c$, and define
$a:=c-h(h+1)$. Then $A(k[h,c],\sigma,a)\cong \U(\mathfrak{sl}_2)$. Under the
obvious isomorphism (choosing $x$, $y$ and $h$ as the standard generators of
$\mathfrak{sl}_2$) the image of the element $c$ corresponds to the Casimir
element.

\item Given $\lambda \in k$, the maximal primitive quotients of
$\U(\mathfrak{sl}_2)$ are the algebras $B_{\lambda}:=\U(\mathfrak{sl}_2)/\langle
c-\lambda\rangle$, cf.~\cite{Dixmier}.  They can also be obtained as
generalized Weyl algebras because $B_{\lambda}\cong A(k[h],\sigma,a=\lambda -
h(h+1))$.  

\end{enumerate}

We will focus on the family of examples $A=A(k[h],a,\sigma)$ with
$a=\sum_{i=0}^na_ih^i\in k[h]$ a non-constant polynomial, and the automorphism $\sigma$
defined by $\sigma(h)=h-h_0$, with $h_0\in k\setminus\{0\}$. 

There is a filtration on $A$ which assigns to the generators $x$ and $y$
degree $n$ and to $h$ degree $2$; the associated graded algebra is, with an
obvious notation, $k[x,y,h]/(yx-a_nh^n)$. This is the coordinate ring of a
Klein surface. We remark that this is a complete intersection, hence it is
a Gorenstein algebra.

There is also a graduation on
$A$, which we will refer to as weight, such that $\deg x=1$, $\deg y=-1$
and $\deg h=0$.

We will denote, for polynomial $a, b\in k[h]$, $\deg a$ the degree of $a$,
$a'=\frac{\partial a}{\partial h}$ the formal derivative of $a$, and
$(a;b)$ the  greatest  common divisor of $a$ and $b$.

Our main results are the following theorems, whose proofs will be given in next
sections.

\begin{teo}\label{thm:hom}
Let $a\in k[h]$ be a non-constant polynomial, $\sigma\in \Aut_k(k[h])$
defined by $\sigma(h)=h-h_0$ with $0\neq h_0\in k$.  Consider
$A=A(k[h],a,\sigma)$ and $n=\deg a$, $d=\deg(a;a')$.
\begin{itemize}
\item  If $(a;a')=1$ (i.e. $d=0$), then $\dim_k\HH_0(A)=n-1$,
$\dim_k\HH_2(A)=1$, and $\HH_i(A)=0$ for $i\neq 0,2$.
\item If $d\geq 1$ then $\dim_k\HH_0(A)=n-1$, $\dim_k\HH_1(A)=d-1$, and
$\dim_k\HH_i(A)=d$ for $i\geq 2$.
\end{itemize}
\end{teo}

\begin{teo}\label{thm:coho}
Let $a\in k[h]$ be a non-constant polynomial, $\sigma\in \Aut_k(k[h])$
defined by $\sigma(h)=h-h_0$ with $0\neq h_0\in k$.  Consider
$A=A(k[h],a,\sigma)$ and $n=\deg a$, $d=\deg(a;a')$.
\begin{itemize}
\item  If $(a;a')=1$ (i.e. $d=0$), then $\dim_k\HH^0(A)=1$,
$\dim_k\HH^2(A)=n-1$, and $\HH^i(A)=0$ for $i\neq 0,2$.
\item If $d\geq 1$ then $\dim_k\HH^0(A)=1$, $\dim_k\HH^1(A)=0$,
$\dim_k\HH^2(A)=n-1$, and $\dim_k\HH^i(A)=d$ for $i\geq 3$.
\end{itemize}
\end{teo}

%%%%%%%%%%%%%%%%%%%%%%%%%%%%%%%%%%%%%%%%%%%%%%%%%%%%%%%%%%
\section{A resolution for $A$ and proof of the first theorem}\label{sect:res}

In this section we will construct a complex of free $A^e$-modules and we
will prove, using an appropriate filtration, that this complex is actually
a resolution of $A$.

The construction of the resolution is performed in two steps.  First we
consider an algebra $B$ above $A$ which has ``one relation less'' than $A$.
Then we use the Koszul resolution of $B$ and obtain a resolution of $A$
mimicking the construction of the resolution for the coordinate ring of an
affine hypersurface done in \cite{burg-vigue}.

Let us consider $V$ the $k$-vector space with basis $\{e_x, e_y, e_h\}$ and
the following complex of free $A^e$-modules:
  \begin{equation}\label{dag}
  0\to A\ot \Lambda^3V\ot A \to A\ot \Lambda^2V\ot A \to
   A\ot V\ot A \to A\ot A \to 0
  \tag{\dag}
  \end{equation}

In order to define the differential, we consider the elements $\lambda_k\in k$ such that
  \[
  \sigma(a)-a=\sum_{k=0}^{n-1}\lambda_kh^k
  \]
and let
  \[
  e_{[x,y]}=\sum_{k,i}\lambda_kh^ie_hh^{k-i-1}, \:
  e_{[x,h]}=-e_x, \:
  e_{[y,h]}=e_y \in A\ot V\ot A.
  \]
For simplicity, in these formulas we have written for example $xe_y$
instead of $x\ot e_y\ot 1\in A\ot V\ot A$ and similarly in the other
degrees.

The  differential in \eqref{dag} is formally, with these notations,  the 
Chevalley-Eilenberg differential; for example,
  \[
  d(\alpha e_x\wedge e_y \beta)
  = \alpha x e_y \beta - \alpha  e_y x \beta
   - \alpha  y e_x \beta + \alpha  e_x y \beta
   - \alpha  e_{[x,y]} \beta.
  \]

\begin{lem}
The homology of the complex \eqref{dag} is isomorphic to $A$ in degree $0$ and
$1$, and zero elsewhere.
\end{lem}

\begin{proof}
We consider the $k$-algebra $B$ 
freely generated by $x,y,h$ modulo relations
  \[
  xh=\sigma(h)x, \qquad
  hy=y\sigma(h), \qquad
  xy-yx=\sigma(a)-a.
  \]
Let $f:=yx-b\in B$ where $b\in k[h]$ is such that $\sigma(b)-b=a$;
observe that $f$ is central in $B$.  This algebra $B$ has been studied by
S.~Smith in \cite{smith}.  Our interest in it comes from the fact that $A$ is
the quotient of $B$ by the two sided ideal generated by $f$. In
particular, $B$ has a filtration induced by the filtration 
on $A$, and it is clear that the associated graded algebra is simply
$k[x,y,h]$.

We claim that a complex similar to \eqref{dag} but with $A$ replaced throughout by $B$
is a resolution of $B$ by free $B^e$-modules.  Indeed, the filtration on
$B$ extends to a filtration on this complex, and the associated graded
object is acyclic, because it coincides with the usual Koszul resolution of
the polynomial algebra $k[x,y,h]$ as a bimodule over itself.

The original complex can be recovered by tensoring over $B$ with $A$ on the
right and on the left, or, equivalently, by tensoring over $B^e$ with
$A^e$.  As a result, the homology of the resulting complex computes
$\Tor_*^{B^e}(B,A^e)\cong \Tor_*^B(A,A)$ (see for example \textsc{IX}.\S4.4 of
\cite{CE}). 

Consider now the free resolution of $A$ as a left $B$-module
  \[
  0\to B\to B\to A\to 0
  \]
with the map $B\to A$ being the natural projection and the other one
multiplication by $f$. 
This can be used to compute $\Tor^B_*(A,A)$, and the proof
of the lemma is finished.
\end{proof}

In order to kill the homology of the complex \eqref{dag}, we consider a
resolution of the following type:
  \begin{equation}
  \xymatrix@C-10pt@-10pt{
  &&0\ar[r]& A\ot \Lambda^3V\ot A \ar[r]& A\ot \Lambda^2V\ot A \ar[r]
  & A\ot V\ot A \ar[r]& A\ot A \ar[r]& 0\\
  &0\ar[r]& A\ot \Lambda^3V\ot A \ar[u]\ar[r]& A\ot \Lambda^2V\ot A\ar[u] \ar[r]& 
  A\ot V\ot A \ar[u]\ar[r]& A\ot A \ar[u]\ar[r]& 0\ar[u]&\\
  0\ar[r]& A\ot \Lambda^3V\ot A \ar[u]\ar[r]& A\ot \Lambda^2V\ot A\ar[u] \ar[r]& 
  A\ot V\ot A \ar[u]\ar[r]& A\ot A \ar[u]\ar[r]& 0\ar[u]&&\\
  &\ar[u] &\ar[u]&\ar[u]&\ar[u] &&
  }
  \tag{\ddag}\label{ddag}
  \end{equation}
The horizontal differentials are the same as before, and the vertical
ones--denoted ``$.df$"--are defined as follows:
  \begin{align*}
  .df:A\ot A &\to A\ot V\ot A\\
  1\ot 1&\mapsto y e_x+e_y x-\sum_{i,k}a_kh^i e_hh^{k-i-1}\\
  \end{align*}
  \begin{align*}
  .df:A\ot V\ot A &\to A\ot \Lambda^2V\ot A\\
  e_x&\mapsto  -e_x\wedge e_y x
  +\sum_{i,k}a_k\sigma(h^i) e_x\wedge e_hh^{k-i-1}\\
  e_y&\mapsto -y e_y\wedge e_x 
  +\sum_{i,k}a_kh^i e_y\wedge e_h\sigma(h^{k-i-1})\\
  e_h&\mapsto -y e_h\wedge e_x-  e_h\wedge e_y x\\
  \end{align*}
  \begin{align*}
  .df:A\ot \Lambda^2V\ot A &\to A\ot \Lambda^3V\ot A\\
  e_y\wedge e_h&\mapsto y e_y\wedge e_h\wedge e_x  \\
  e_x\wedge e_h&\mapsto  e_x\wedge e_h\wedge e_y  x\\
  e_x\wedge e_y&\mapsto 
  -\sum_{i,k}a_k\sigma(h^i) e_x\wedge e_y\wedge e_h  \sigma(h^{k-i-1})
  \end{align*}
  
\begin{prop}
The total complex associated to \eqref{ddag} is a resolution of $A$ as
$A^e$-module.
\end{prop}

\begin{proof}
That it is a double complex follows from a straightforward computation. 
We consider again the filtration on $A$ and the filtration induced by it 
on \eqref{ddag}. All maps 
respect it, so it will suffice to see that the associated graded complex 
is a resolution of $\gr A$.
Filtering this new complex by rows, we know that the homology of the rows
compute $\Tor_*^{\gr B}(\gr(B),\gr(B))$. The only thing to be checked now
 is that the differential on the
$E^1$ term can be identified to $.df$, and this is easily done.
\end{proof}

In order to compute $\HH_*(A)$ we can compute the homology of the complex
$(A\ot_{A^e}(A\ot \Lambda^*V\ot A),.df,d_{CE})\cong (A\ot
\Lambda^*V,.df,d_{CE})$. This is a double complex which can be filtered by the
rows, as usual, so we obtain a spectral sequence converging to the homology of the total
complex. Of course, the first term is just the homology of the rows.

The computation of Hochschild homology can be done in a direct way; however, it
is worth noticing that this procedure can be considerably reduced. Let
$X_*=(A\ot\Lambda^*V,d)$ be the complex obtained by tensoring the rows in
\eqref{ddag} with $A$ over $A^e$, and let $X^0_*$ be the zero weight
component of $X_*$.

\begin{prop}\label{prop:reduction}
The inclusion map $X^0_*\to X_*$ induces an isomorphism in homology.
\end{prop}

\begin{proof}
Let us define the map $s:A\ot \Lambda^*V\to A\ot \Lambda^{*+1}V$ 
by $s(w\ot v_1\wedge \dots\wedge v_k):= w\ot v_1\wedge \dots\wedge v_k\wedge
e_h$. A computation shows that
  \begin{equation}\label{eq:htpy}
  (d_{CE}s+sd_{CE})(w\ot v_1\wedge \dots\wedge v_k)=
  weight(w\ot v_1\wedge \dots\wedge v_k)
  .(w\ot v_1\wedge \dots\wedge v_k).
  \end{equation}
Since $\chr(k)=0$, this ``Euler'' map is an isomorphism for non-zero
weights, but it is the zero map in homology, because \eqref{eq:htpy} shows that it
is homotopic to zero.
\end{proof}

%%%%%%%%%%%%%%%%%%%%%%%%%%%%%%%%%%%%%%%%%%%%%%%%
\subsection{The term $E^1$}

%%%%%%%%%%%%%%%%%%%%%%%%%%%%%%%%%%%%%%
\paragraph{Computation of $\HH_0(A)$.}
We remark that in the graduation by weight $A=\oplus_{n\in\ZZ}A_n$ we have 
$A_0=k[h]$, and, for $n>0$, $A_n=k[h]x^n$ and $A_{-n}=k[h]y^n$.

We have to compute $\HH_0(A)=A/[A,A]=A/([A,x]+[A,y]+[A,h])$, and this is, according to
proposition \ref{prop:reduction}, the same as $A_0/([A_{-1},x]+[A_{1},y]+[A_0,h])$.

Since $A_0=k[h]$, $[A_0,h]=0$.  Because $A_{-1}=k[h]y$, a
system of linear generators of $[A_{-1},x]$ is given by commutators of the form
$[h^iy,x]=h^ia-\sigma(h^ia)=(\Id-\sigma)(h^ia)$.
On the other hand, $A_1=k[h]x$, so $[A_1,y]$ is spanned by the
$[h^jx,y]=h^j\sigma(a)-\sigma^{-1}(h^j)a=(\Id-\sigma)(-\sigma^{-1}(h^j )a)$
for $j\geq0$.  As a consequence, $[A_1,y]+[A_{-1},x]=[A_{-1},x]$ is the
subspace of $A_0=k[h]$ of all polynomials $pa-\sigma(pa)$ with $p\in k[h]$.

The $k$-linear map $\Id-\sigma:k[h]\to k[h]$ is an epimorphism, and its kernel
is the one dimensional subspace consisting of constant polynomials. The subspace
of multiples of $a$ has codimension $\deg a=n$, so the image of the restriction
of $\Id-\sigma$ to this subspace has codimension $n-1$.  Then we conclude that
$\dim_k\HH_0(A)=n-1$; a basis is given for example by the set of homology
classes $\{[1], [h], [h^2],\dots ,[h^{n-2}]\}$.

%%%%%%%%%%%%%%%%%%%%%%%%%%%%%%%%%%%%%%
\paragraph{Homology of the row in degree $1$.}
Recall that we only have to consider the subcomplex of elements of weight zero.
Let us then suppose that $c=sye_x+txe_y+ue_h$ is a $1$-cycle in the row complex
of weight zero.  This implies that
  \[
  d(sye_x+txe_y+ue_h)=
  	\sigma\left( (\sigma^{-1}(t)-s)a \right) - (\sigma^{-1}(t)-s)a = 0.
  \]
As a consequence, $(\sigma^{-1}(t)-s)a\in \Ker(\Id-\sigma)=k$. But
$a$ is not a constant polynomial, so $\sigma^{-1}(t)-s=0$. In other words,
$s=\sigma^{-1}(t)$ and the cycle can be written in the form
  \[
  \sigma^{-1}(t)ye_x+txe_y+ue_h.
  \]
The horizontal boundary of a $2$-chain $pe_x\wedge e_y +qye_x\wedge e_h +rxe_y\wedge e_h$ is
  \begin{multline*}
  d_{CE}(pe_x\wedge e_y +qye_x\wedge e_h +rxe_y\wedge e_h)= \\
  =(p-\sigma(p))xe_y + (\sigma^{-1}(p)-p)ye_x
  +\left(-p (\sigma(a')-a')
  +(q-\sigma^{-1}(r))a - \sigma((q-\sigma^{-1}(r))a)\right)e_h.
  \end{multline*}
We can choose $p$ such that $\sigma(p)-p=t$, so that, adding $d(pe_x\wedge e_y)$
to $c$, we obtain cycle homologous to $c$, in which the only eventually non-zero
coefficient is the one corresponding to $e_h$.  We can then simply assume that
$c$ is of the form $ue_h$ to begin with, and we want to know if it is a boundary
or not.

The equation $d(pe_x\wedge e_y +qye_x\wedge e_h +rxe_y\wedge e_h)=ue_h$
implies $p=\sigma(p)$, so $p\in k$, and
  \begin{equation}\label{yyy}
  u+ \left(-(\sigma(a')-a')p \right)
  =- \sigma((q-\sigma^{-1}(r))a) +(q-\sigma^{-1}(r))a.
  \end{equation}
If $n=\deg a=1$, then we are in the special case of the usual Weyl algebra.
In this case $a'\in k$ and $\sigma(a')-a'=0$, so there is one term less in the
left hand side of \eqref{yyy} and the homology of the row in this degree is zero.

Suppose now $n\geq 2$;  if $p=0$ then $u\in\Im((\Id-\sigma)|_{a.k[h]})$, and
we have, as in degree zero, that $\{[1],\dots,[h^{n-2}]\}$ is a basis of the
quotient.  If $p\neq 0$, we have to mod out a $(n-1)$-dimensional space by the
space spanned by a non-zero element, so we obtain a $(n-2)$-dimensional space.  We
notice that $\deg(\sigma(a')-a')=\deg a-2$, and since $n\geq 2$, the element
$\sigma(a')-a'$ is a non-zero element of $k[h]/(\Id-\sigma)(a.k[h])$.

%%%%%%%%%%%%%%%%%%%%%%%%%%%%%%%%%%%%%%
\paragraph{Homology of the row in degree $2$.}
The boundary of a weight zero element in degree two\\
$w:=se_x\wedge e_y + tye_x\wedge e_h +uxe_y\wedge e_h$
is:
  \[
  d(w)=
  (s-\sigma(s))x e_y -
  (s-\sigma^
{-1}(s))y e_x
   +\left( (a'-\sigma(a'))s
  +(t-\sigma^{-1}(u))a -
  \sigma((t-\sigma^{-1}(u))a) \right)e_h.
  \]
If $d(w)=0$, we must have that $s=\sigma(s)$, so $s\in k$, and that
  \[
  s(\sigma(a')-a')=-(t-\sigma^{-1}(u))a -\sigma((t-\sigma^{-1}(u))a).
  \]
If $s\neq 0$, then the expression on the left is a polynomial of degree $n-2$,
and the degree of the polynomial on the right is (if it is not the zero
polynomial) $n+\deg(t-\sigma^{-1}(u))-1$.  This is only possible if both sides are
zero and we see that $\sigma(t)=u$.

We mention that in the case $\deg a=1$ (i.e. the usual Weyl algebra), the
expression on the left is always zero independently of $s$, so the argument is
not really different in this case.

Now we compute the $2$-boundaries: as $p$ varies in $k[h]$, they are the
elements
  \begin{align*}
  d(p e_x \wedge e_y \wedge e_h)&=
  [p,x]e_y\wedge e_h -[p,y]e_x\wedge e_h +[p,h]e_x\wedge e_y\\
  &=(p-\sigma(p))x e_y\wedge e_h - (p-\sigma^{-1}(p))y e_x\wedge e_h.
  \end{align*}
Given $u$, there is a $p$ such that $(\Id-\sigma)(p)=u$,
and this $p$ automatically satisfies $\sigma^{-1}(p)-p=t$.

We remark that the coefficient corresponding to
$e_x\wedge e_y$ in a $0$-weight boundary, is always zero.
As a consequence, in the case of the usual Weyl algebra,
the class of $e_x\wedge e_y$ is a generator of the homology.
On the other hand, if $n\geq 2$ the homology is zero.

%%%%%%%%%%%%%%%%%%%%%%%%%%%%%%%%%%%%%%
\paragraph{Homology of the row in degree $3$.}
The homology in degree three is the kernel of the map
$A\ot \Lambda^3V \to A\ot \Lambda^2V$ given by
  \[
  w\mapsto [w,x]e_y\wedge e_h -[w,y]e_x\wedge e_h +[w,h]e_x\wedge e_y.
  \]
It is clearly isomorphic to the center of $A$, which is known to be $k$ (see
for example \cite{Bav}). A basis of the homology is given by the class of
$e_x\wedge e_y\wedge e_h$.

%%%%%%%%%%%%%%%%%%%%%%%%%%%%%%%%%%%%%%
\paragraph{Summary.} 
We summarize the previous computations in the following 
table showing the dimensions of the vector spaces in the term $E^1$. In each
case, the boxed entry has coordinates $(0,0)$.
\[
\begin{array}{ccc}

\begin{array}{ccccccc}
 & &   &1    &0  & n-2&\fbox{$n-1$}\\
 & &1  &0    &n-2&n-1&\\
 &1& 0 & n-2 &n-1&   &\\
1&0&n-2& n-1 &   &   &\\
\end{array}
&
\qquad
&
\begin{array}{ccccccc}
 & &   &1    &1  & 0&\fbox{$0$}\\
 & &1  &1    &0&0&\\
 &1& 1 & 0 &0&   &\\
1&1&0& 0 &   &   &\\
\end{array}
\\
\\
n\geq2 && \text{The Weyl algebra ($n=1$)}
\end{array}
\]

%%%%%%%%%%%%%%%%%%%%%%%%%%%%%%%%%%%%%%%%%%%%%%%%
\subsection{The term $E^2$}

The differential $d^1$ corresponds to the vertical differential in the original
complex. Let $n\geq 2$. The only relevant component is the map $.df:A\to A\ot V$;
we recall that it is defined by
  \[
  .df(b)= bye_x +\sigma(b)xe_y-ba'e_h.
  \]
Adding $d_{CE}(pe_x\wedge e_y)$, where $p$ is such that 
$b=\sigma(p)-p$, we see that the expression $bye_x +\sigma(b)xe_y-ba'e_h$ is
homologous to $\left(-\sigma\left(\sigma^{-1}(p) a'\right)+\left(\sigma^{-1}(p)
a'\right)\right)e_h$, and, since the homology of the row in the place
corresponding to $A\ot V$ is isomorphic to $k[h]/(\Id-\sigma)(a.k[h])e_h$, we
conclude that the cokernel of the first differential of the spectral sequence
(in the same place) is isomorphic to $k[h]/(\Id-\sigma)(a.k[h] +a'k[h])e_h$.
The subspace $ak[h]+a'k[h]$ has codimension $d=\deg(a;a')$, so
$(\Id-\sigma)(ak[h]+a'k[h])$ has codimension $d-1$ (or zero if $d=0$). 

By linear algebra arguments, the dimension of the kernel of this
differential is $d$. The corresponding table at this step of the
spectral sequence is the following:
  \[
  \begin{array}{ccc}
    \begin{array}{ccccccc}
     & &   &1    &0  & d-1&\fbox{$n-1$}\\
     & &1  &0    &d-1&d&\\
     &1& 0 &d-1  &d&   &\\
    1&0&d-1&d  &   &   &
    \end{array}
  &
  \qquad
  &
    \begin{array}{ccccccc}
     & &   &1    &0  & 0&\fbox{$n-1$}\\
     & &1  &0    &0&1&\\
     &1& 0 &0  &1&   &\\
    1&0&0&1  &   &   &
    \end{array}
  \\
  &&\\
  d\geq1 && d=0
  \end{array}
  \]

In case $n=1$, the homology of the Weyl algebra is well-known (see for example
\cite{kassel} or \cite{sri}), but for completeness we include
it. The only relevant differential is the one corresponding to the map $A\ot
\Lambda^2V\to A\ot \Lambda^3V$. The generator of the homology of the row in the
place corresponding to $A\ot\Lambda^2V$ is $e_x\wedge e_y$, and $.df(e_x\wedge
e_y)=-\sigma(a')e_x\wedge e_y\wedge e_h$. But $\deg a=1$ so $a'$ is a non-zero
constant; as a consequence $d_1$ is an epimorphism and hence an isomorphism. The table of the
dimensions is in this case:
  \[
  \begin{array}{c}
  \begin{array}{ccccccc}
   & &   &0    &1  & 0&\fbox{$0$}\\
   & &0  &0    &0&0&\\
   &0& 0 &0  &0&   &\\
  0&0&0&0  &   &   &\\
  \end{array}
  \\
  \\
  n=1
  \end{array}
  \]
We recover thus known results.

%%%%%%%%%%%%%%%%%%%%%%%%%%%%%%%%%%%%%%%%%%%%%%%%
\subsection{The term $E^3$}\label{subsect:e3}

Since $d^2$ has bidegree $(-2,1)$, its only eventually non-zero
component has as target a vector space of dimension one, and there are
two possibilities: either it is zero, or it is an epimorphism.

In order to decide whether it is an epimorphism, it is sufficient to determine
if the element $e_x\wedge e_y\wedge e_h$ is a coboundary or not.
In other words, we want to know if there exist
$z_1=\alpha xe_y\wedge e_h+\beta ye_x\wedge e_h+\gamma e_x\wedge e_y$
and $z_2=p$ such that $z_1.e_f=e_x\wedge e_y \wedge e_h$ and  
  \begin{equation}\label{eq:cond}
  d_{CE}(z_1)+.df(z_2)=0.
  \end{equation}
We have that $.df(z_1)=(((\alpha-\sigma(\beta))\sigma(a)-\gamma\sigma( a'))e_x\wedge e_y\wedge e_h$;
so $.df(z_1)=e_x\wedge e_y\wedge e_h$ if and only if
  \begin{equation}\label{eq:natural}
  (\alpha-\sigma(\beta))\sigma(a)-\gamma\sigma( a')=1,
  \end{equation}
if and only if
  \[
  (\sigma^{-1}(\alpha)-\beta)a-\sigma^{-1}(\gamma) a'=1.
  \]
A necessary and sufficient condition for a solution to this equation to exist
is that $(a;a')=1$, in other words, that $a$ have only simple roots.

If this is the case, let $(\alpha,\beta,\gamma)$ be a solution.
We have
  \begin{multline*}
  d_{CE}(z_1)=
  (\gamma- \sigma(\gamma))xe_y- (\gamma- \sigma^{-1}(\gamma))ye_x+\\
  +\left(-\gamma(\sigma(a')-a')+
  ( \beta- \sigma^{-1}(\alpha))a-
  \sigma(( \beta- \sigma^{-1}(\alpha))a)\right)e_h;\qquad\qquad\qquad
  \end{multline*}
and using  \eqref{eq:natural} we see that this is equal to
  \[
  (\gamma- \sigma(\gamma))xe_y- (\gamma- \sigma^{-1}(\gamma))ye_x
  +(\gamma-\sigma^{-1}(\gamma))a'e_h.
  \]
On the other hand, $.df(z_2)=.df(p)= pye_x +\sigma(p)xe_y-pa'e_h$. It
is then enough to choose $p=\sigma^{-1}(\gamma)-\gamma$ to have equation
\eqref{eq:cond} satisfied.

We conclude that $d_2$ is an epimorphism if $(a;a')=1$, and zero if not. 

We can summarize the results of the above computations
in the following table containing the dimensions
of $\HH_p(A)$:
  \[
  \begin{array}{||c|c|c||}
  \hline
  p&(a:a')=1&\deg((a:a'))=d\geq 1\\ 
  \hline
  \hline
  0           &  n-1   &         n-1  \\
  \hline
  1           & 0      &         d-1  \\      
  \hline
  2           & 1      &         d    \\     
  \hline
  \geq 3      &0       &          d   \\       
  \hline
  \end{array}
  \]
We note that this proves theorem \ref{thm:hom}.

%%%%%%%%%%%%%%%%%%%%%%%%%%%%%%%%%%%%%%%%%%%%%%%%%%%%%%%%%%
\section{Cohomology}\label{sect:coho}

The aim of this section is to compute the Hochschild cohomology of GWA. Once we
have done this, we compare the obtained dimensions with duality results.

We use resolution \eqref{ddag} of $A$ as an $A^e$-module to compute
cohomology.  We apply the functor $\Hom_{A^e}(-,A)$ and make the
following identifications:
  \[
  \Hom_{A^e}(A\ot \Lambda^kV\ot A,A)\cong
  \Hom(\Lambda^kV,A)\cong
  (\Lambda^kV)^*\ot A\cong
  \Lambda^{3-k}V\ot A.
  \]
Here we are identifing $(\Lambda^{k}V)^*\cong \Lambda^{3-k}V$ using the
pairing $\Lambda^{k}V\ot\Lambda^{3-k}V\to\Lambda^3V\cong k$ given by
exterior multiplication.
Using superscripts for the dual basis, the correspondence between the basis
of $(\Lambda^{3-k}V)^*$ and the basis of $\Lambda^k V$ is:
\begin{align*}
e^h&\mapsto e_x\wedge e_y   
  & e^x\wedge e^h&\mapsto-e_y 
  & 1&\mapsto e_x\wedge e_y\wedge e_h \\
e^y&\mapsto-e_x\wedge e_h 
  &e^x\wedge e^y&\mapsto e_h
  \\
e^x&\mapsto e_y\wedge e_h
    &e^y\wedge e^h&\mapsto e_x
\end{align*}

In this way, we obtain the following double complex, whose total homology
computes $\HH^*(A)$:
  \[
  \xymatrix@C-10pt@R-10pt{
  0\ar[r]& A\ot \Lambda^3V \ar[d]\ar[r]& A\ot \Lambda^2V\ar[d] \ar[r]& A\ot V \ar[d]\ar[r]& A \ar[d]\ar[r]& 0\ar[d]&&\\
  &0\ar[r]& A\ot \Lambda^3V \ar[d]\ar[r]& A\ot \Lambda^2V\ar[d] \ar[r]& A\ot V \ar[d]\ar[r]& A \ar[d]\ar[r]& 0\ar[d]&\\
  &&0\ar[r]& A\ot \Lambda^3V \ar[d]\ar[r]& A\ot \Lambda^2V\ar[d] \ar[r]& A\ot V \ar[d]\ar[r]& A \ar[d]\ar[r]& 0\\
  &&&&&&&}
  \]
The horizontal differentials are, up to sign, {\em exactly} the same as in
homology. The vertical differentials are also essentially the same; for
example, $.df^*:A\to A\ot V$ is such that
  \[
  .df^*(b)=yb\ot e_x+
  bx\ot e_y-\sum_{i,k}h^{k-i-1}ba_kh^i\ot e_h.
  \]

\begin{remark}
The differences (and similarities) between the above formulas and the
corresponding ones in homology may be explained as follows.  Given the map
$A^e\to A^e$ defined by $a\ot b\mapsto az\ot wb$, the induced maps on the
tensor product and $\Hom$ are related in the following way: when we use the
tensor product functor we obtain:
  \begin{align*}
  A\cong A\ot_{A^e}A^e&\to  A\ot_{A^e}A^e \cong A\\
  b&\mapsto  wbz.
  \end{align*}
When, on the other hand, we use the $\Hom$ functor we get:
  \begin{align*}
  A\cong \Hom_{A^e}(A^e,A)&\to \Hom_{A^e}(A^e,A) \cong A\\
  b&\mapsto  zbw.
  \end{align*}
\end{remark}

This fact implies that we already know the homology of the rows, up to
reindexing. However it is worth noticing that there is a change of degree
with respect to the previous computation (now degree increases from left to
right).  Schematically, the dimensions of these homologies are:
  \[
  \begin{array}{ccccccc}
  \fbox{$1$}    &0  & n-2& n-1&&&\\
  &1    &0  & n-2& n-1&&\\
  &&1    &0  & n-2& n-1&\\
  &&&1    &0  & n-2& n-1
  \end{array}
  \]
From this, it follows that $\dim_k\HH^0(A)=1$ and $\dim_k\HH^1(A)=0$,
independently of the polynomial $a$. Also $\dim_k\HH^2(A)=n-1=\deg a-1$,
because of the form of the $E_1$ term in the spectral sequence.

As before, there are two different cases: either \textsc{(i)} $(a;a')=1$ or 
\textsc{(ii)} $1\leq \deg(a;a')=d\leq n-1$. The following tables give the
dimensions of the spaces in the $E_2$ terms of the spectral sequence, in
both situations:
  \[
  \begin{array}{ccc}
  \begin{array}{ccccccc}
  \fbox{$1$}    &0  & n-2& 1&&&\\
  &1    &0  & 0& 1&&\\
  &&1    &0  & 0& 1&\\
  &&&1    &0  &0 & 1\\
  \end{array}
  &\qquad&
  \begin{array}{ccccccc}
  \fbox{$1$}    &0  & n-2& d&&&\\
  &1    &0  & d-1& d&&\\
  &&1    &0  & d-1& d&\\
  &&&1    &0  &d-1 &d \\
  \end{array}
  \\
  \\
  \text{\textsc{(i)}}&&\text{\textsc{(ii)}}
  \end{array}
  \]
In each case, the differential $d_2$ is the same, up to our
identifications, as the one considered in section \ref{subsect:e3}. As a
consequence, we have: in case \textsc{(i)}, the $E_3$ term has the form
  \[
  \begin{array}{ccccccc}
  \fbox{$1$}    &0  & n-2& 0&&&\\
  &1    &0  & 0& 0&&\\
  &&0    &0  & 0& 0&\\
  &&&0    &0  &0 & 0\\
  \end{array}
  \]
so $E_\infty=E_3$; in case \textsc{(ii)}, $d_2=0$, and we see that
$E_\infty=E_2$.

We summarize the results in the following table containing
the dimensions of $\HH^p(A)$:
  \[
  \begin{array}{||c|c|c||}
  \hline
  p&(a;a')=1&\deg(a;a')=d\geq 1\\
  \hline
  \hline
  0           &  1   &      1  \\
  \hline
  1           & 0      &     0  \\
  \hline
  2           & n-1      &    n-1   \\
  \hline
  \geq 3      &0       &          d   \\
  \hline
  \end{array}
  \]
This proves theorem \ref{thm:coho} concerning cohomology.

%\rem encontrar y poner  $n-1$ deformaciones de $A$. Observar
%que (salvo el caso Weyl) ninguna es rigida. Todas las
%regulares (Weyl incluida) tienen $H^3=0$.

It is clear that, when the polynomial $a$ has multiple roots, there is no
duality between Hochschild homology and cohomology, contrary to what one might
expect after \cite{VdB}.  This is explained by the fact that in this case the
algebra $A^e$ has infinite left global dimension; in this situation, theorem 1
in \cite{VdB} fails: one cannot in general replace $\ot$ by $\ot^L$ in the first
line of Van den Bergh's proof. One can retain, however, the conclusion in the
theorem if one adds the hypothesis that either the $A^e$-module $A$ or the
module of coefficients has finite projective dimension.
This is explained in detail in \cite{VdBerratum}.

%%%%%%%%%%%%%%%%%%%%%%%%%%%%%%%%%%%%%%%%%%%%%%%%%%%%%%%%%%
\section{Invariants under finite group actions}\label{sect:inv}

The algebraic torus $k^*=k\setminus\{0\}$ acts on generalized Weyl algebras
by diagonal automorphisms. More precisely, given $w\in k^*$,
there is an automorphism of algebras uniquely determined by
  \[
  x\mapsto wx, \qquad y\mapsto w^{-1}y, \qquad h\mapsto h.
  \]
This defines a morphism of groups $k^*\hookrightarrow\Aut_k(A)$ which we will
consider as an inclusion. The automorphism defined by $w\in k^*$ is of finite order if
and only if $w$ is a root of unity, and, in this case, the subalgebra of
invariants can easily be seen to be generated by $\{h,x^m,y^m\}$, where $m$ is
the order of $w$.

The following lemma is a statement of the fact that the process of taking
invariants with respect to finite subgroups of $k^*$ for this action does
not lead out of the class of GWA.  This enables us to obtain almost
immediately the Hochschild homology and cohomology of the invariants.

\begin{lem}\label{lemma:inv}
Let $A=A(k[h],\sigma, a)$ be a generalized
Weyl algebra and let $G:=\ZZ/r.\ZZ$ act on $A$ by powers of
the diagonal automorphism induced by a primitive $r$-th root
of unity. The subalgebra of invariants $A^G$ is isomorphic
to the generalized Weyl algebra $A^G=A(k[H],\tau,\wt{a})$,
where $\tau(H)=H-1$
and $\wt{a}(H)=\sigma^{-r+1}(a)(rH)\cdots\sigma^{-1}(a)(rH)a(rH)$.
\end{lem}

\begin{proof}
We know that $A^G=\langle h, x^r,y^r \rangle$.  Let us write $X:=x^r$,
$Y:=y^r$ and $H:=h/r$. Then $XH=x^rh/r=\sigma^r(h/r)x^r=\tau(H)X$ and similarly
$HY=Y\tau(H)$. Now
  \[
  YX=y^rx^r=y^{r-1}yxx^{r-1}=y^{r-1}a(h)x^{r-1}= \sigma^{-r+1}(a)(h)y^{r-1}x^{r-1}
	= \sigma^{-r+1}(a)(rH)y^{r-1}x^{r-1}.
  \]
so clearly the equality $y^rx^r=\wt a(H)$ follows by induction on $r$.
\end{proof}

The idea to compute (co)homology of $A^G$ is to replace it by the crossed
product $A*G$. This change does not affect the homology provided that $A^G$
and $A*G$ are Morita equivalent; this is discussed in detail in
\cite{afls}. In particular, this is the case when the polynomial $a$ has 
no pair of different roots conjugated to each 
other by $\sigma$---that is, there do not
exist $\mu\in\CC$ and $j\in\ZZ$ such that $a(\mu)=a(\mu+j)=0$---because in this
situation the algebra $A=A(k[h],\sigma,a)$ is simple, as proved by Bavula
in~\cite{Bav}. We state this as

\begin{prop}\label{prop:desc}
Let $a\in k[h]$ be a polynomial such that no pair of its roots are
conjugated by $\sigma$ in the sense explained above.
Let $G$ be any finite subgroup of $\Aut_k(k[h])$. Then there are
isomorphisms
  \[
  \HH^*(A^G)
  \cong H^*(A,A*G)^G
  \cong \bigoplus_{\cl{g}\in\cl{G}}H^*(A,Ag)^{\Z_g},
  \]
where the sum is over the set $\cl G$ of conjugacy classes $\cl g$ of $G$,
and, for each $g\in G$, $\Z_g$ is the centralizer of $g$ in $G$.
Also, there are duality isomorphisms $\HH_*(A^G)\cong\HH^{2-*}(A^G)$.
\end{prop}

\begin{proof}
Given the hypotheses in the statement, we are in a situation similar
to the one considered in \cite{afls}. The proposition follows from the
arguments presented there.
The last part concerning homology follows from the duality theorem of
\cite{VdB}, since the global dimensional of $A^G$ is finite; see also section~7
in~\cite{afls}.
\end{proof}

Under appropriate conditions on the $g\in G$, we are able to
compute the $\Z_g$-module $H_*(A,Ag)$. This module has always
finite dimension as $k$-vector space (see proposition \ref{prop:coef})
and the action of $\Z_g$ is determined by an element $\Om\in \Aut_k(A)$
(see proposition \ref{prop:action}). This automorphism $\Om$ is a generalization
of the Cartan involution of ${\frak sl}_2$, and is explained in more detail
in section \ref{section:cartan}. We state the theorem, but since we need to know 
some facts about
the group $\Aut_k(A)$, its proof will finish in
section \ref{section:end}.

\begin{teo}\label{teoG} 
Let us consider a GWA $A=A(k[h],\sigma, a)$ which is simple.
Let $G\subset \Aut_k(A)$ be a finite subgroup such that every
element $g$ of $G$ is conjugated in $\Aut_k(A)$ to an element
in the torus $k^*$. 
Let us define 
$a_1:=\#\{\cl{g}\in\cl{G}\setminus\{\Id\}\hbox{ such that } \Om\notin\Z_g\}$ and
$a_2:=\#\{\cl{g}\in\cl{G}\setminus\{\Id\}\hbox{ such that }\Om\in\Z_g\}$. 
We have that
    \[
    \dim_k\HH^p(A^G)=
    \begin{cases}
    1 & \text{if $p=0$}\\
    (n-1)+na_1+[(n+1)/2]a_2 & \text{if $p=2$}\\
    0 & \text{if $p=1$ or $p>2$}\\
    \end{cases}
    \]
\end{teo} 
\begin{remark} In particular, if the action of $\Z_g$ is trivial, the
formula for $\HH^2$ gives $\dim \HH^2(A^G)=n.\#\cl{G}-1$.
\end{remark}

\begin{proof}
Using the hypotheses and the above proposition, the proof
will follow from the computation of the dimensions of $H^*(A,Ag)$ (proposition
\ref{prop:coef}) and the characterization of the action
(proposition \ref{prop:action}).
\end{proof}

We state the following proposition for automorphisms $g$ of
$A$ diagonalizable but not necessarily of finite order, 
although we do not need such generality.

\begin{prop}\label{prop:coef}
Let $g\in \Aut_k(A)$ different from the
identity and conjugated to an element of $k^*$. Then $H^0(A,Ag)=H^1(A,Ag)=0$,
$\dim_kH^2(A,Ag)=n$, and $\dim_kH^*(A,Ag)=d$ for each $*>2$, where $d=\deg(a;a')$.
Also, $\dim_kH_0(A,Ag)=\deg a=n$, and, for all $*>0$,
$\dim_kH_*(A,Ag)=d$. 
\end{prop}

We can assume that $g\in\Aut_k(A)$ is in fact in $k^*$, since, by
Morita invariance, $H^*(A,Ag)\cong H^*(A,Ahgh^{-1})$ for all
$h\in\Aut_k(A)$.

The groups $H^*(A,Ag)$ can be computed using the complex obtained by
applying to the resolution \eqref{ddag} the functor $\Hom_{A^e}({-},Ag)$;
it can be identified, using the same idea as 
in section \ref{sect:coho}, to the double complex
  \[
  \xymatrix@C-10pt@R-10pt{
  0\ar[r]& Ag\ot \Lambda^3V \ar[d]\ar[r]& Ag\ot \Lambda^2V\ar[d] \ar[r]& Ag\ot V \ar[d]\ar[r]& Ag \ar[d]\ar[r]& 0\ar[d]&&\\
  &0\ar[r]& Ag\ot \Lambda^3V \ar[d]\ar[r]& Ag\ot \Lambda^2V\ar[d] \ar[r]& Ag\ot V \ar[d]\ar[r]& Ag \ar[d]\ar[r]& 0\ar[d]&\\
  &&0\ar[r]& Ag\ot \Lambda^3V \ar[d]\ar[r]& Ag\ot \Lambda^2V\ar[d] \ar[r]& Ag\ot V \ar[d]\ar[r]& Ag \ar[d]\ar[r]& 0\\
  &&&&&&&}
  \]
This complex is graded (setting $weight(g)=0$) in an analogous 
way to the complex which computes $\HH^*(A)$. It is straightforward to verify 
that the homotopy defined in \ref{prop:reduction} may be also used in this case. 
As a consequence, the cohomology of the rows is concentrated in weight zero.

\subsection{The term $E_1$}

\paragraph{Computation of $H^0(A,Ag)$.}

Let $p\in k[h]$ and assume $pge_x\wedge e_y\wedge e_h\in \Ker(d_{CE})$, that
is, that
  \[
  (\sigma(p)-wp)xge_y\wedge e_h-
  (\sigma^{-1}(p)-w^{-1}p)yge_x\wedge e_h+
  0e_x\wedge e_y=0
  \]
Then $(\sigma-w.\Id)(p)=0$, so $p=0$, and the cohomology in degree zero
vanishes.

%%%%%%%%%%%%%%%%%%%%%%%%%%%%%%%%%%%%%
\paragraph{Homology of the rows, degree $1$.}

Given $u$, $v$, $t$ in  $k[h]$, a computation shows that
$uge_x\wedge e_y+ vyge_x\wedge e_h+ txge_y\wedge e_h\in
\Ker(d_{CE})$, if and only if
  \begin{multline*}
  (\sigma^{-1}(u)-w^{-1}u)yge_x +(wu-\sigma(u))xge_y+\\
  +(w^{-1}t\sigma(a)-\sigma^{-1}(t)a -u(\sigma(a')-a') +wva-\sigma(va))g
  e_h=0.\qquad\qquad\qquad
  \end{multline*}
Using again that $\sigma-w.\Id$ is an isomorphism, we conclude that
in order for the coefficient of $e_y$ to vanish, $u$ must be zero.
Looking at the coefficient of $e_h$, we obtain the only other condition, it
is
  \[
  w^{-1}t\sigma(a)-\sigma^{-1}(t)a +wva-\sigma(va)=
  (\sigma-w.\Id)((w^{-1}\sigma^{-1}(t)-v)a)=0,
  \]
so $w^{-1}\sigma^{-1}(t)-v=0$. 

We conclude that any $1$-cocycle is of the form
  \[
  d_{CE}(pge_x\wedge e_y\wedge e_h)=
  vyge_x\wedge e_h+ w\sigma(v)xge_y\wedge e_h,
  \]
where $p\in k[h]$ is chosen so that $\sigma(p)-wp=w\sigma(v)$.
It follows immediately that the cohomology of the rows in degree $1$ is zero.

%%%%%%%%%%%%%%%%%%%%%%%%%%%%%%%%%%%%%
\paragraph{Homology of the rows, degree $2$.}

The $2$-coundaries are expressions of the form
  \begin{multline}\label{eq:2bd}
  (\sigma^{-1}(u)-w^{-1}u)yge_x +(wu-\sigma(u))xge_y+ \\
  +(w^{-1}t\sigma(a)-\sigma^{-1}(t)a -u(\sigma(a')-a') +wva-\sigma(va))g e_h
  \end{multline}
A $2$-cochain $\alpha=pyge_x+qxge_y+rge_h$, with $p$, $q$, $r\in k[h]$ is a
cocycle if and only if $q=w\sigma(p)$ holds; in particular, this imposes no
conditions on $r$.  We can then assume that $p=q=0$ in $\alpha$ 
because one can add to $\alpha$ a coboundary of the form $d(ug e_x\wedge e_y)$ with $u\in
k[h]$. We want now to decide when such a $2$-cocycle is a
coboundary. In view of \eqref{eq:2bd} and the fact that
$\sigma-w.\Id$ is an isomorphism, we see at once $u$ must be zero.
We are reduced to solve the equation
  \[
  w^{-1}t\sigma(a)-\sigma^{-1}(t)a 
  +wva-\sigma(va)=(\sigma-w.\Id)((w^{-1}\sigma^{-1}(t)-v)a)= r.
  \]
This can be solved if and only if
$(\sigma-w.\Id)^{-1}(r)$ is a multiple of $a$, so the codimension
of the subspace of solutions is $\deg a=n$, in other words,
the dimension of the cohomology of the rows in degree $2$ is $n$.

%%%%%%%%%%%%%%%%%%%%%%%%%%%%%%%%%%%%%
\paragraph{Homology of the rows, degree $3$.}

The coboundaries of weight zero are of the form
  \[
  d(pyge_x+qxge_y+rge_h)=
  (wpa-\sigma(pa)+w^{-1}q\sigma(a)-\sigma^{-1}(q)a)g=
  (\sigma-w.\Id)((w^{-1}\sigma^{-1}(q)-p)a)g
  \]
with $p$, $q$, $r$ in $k[h]$. Every polynomial in $k[h]$ can be written as
$w^{-1}\sigma^{-1}(q)-p$ for some $p, q\in k[h]$, so, since the map
$\sigma-w.\Id:k[h]\to k[h]$ is an isomorphism, we see that the dimension of the
cohomology of the row complex in degree $3$ is equal to $\dim_kk[h]/ak[h]=\deg
a=n$.

%%%%%%%%%%%%%%%%%%%%%%%%%%%%%%%%%%%%%
\subsection{The term $E_2$}

In view of the above computations, the dimensions of the components of the 
$E_1$-term of the spectral sequence are as follows:
  \[
  \begin{array}{cccccc}
  \fbox{$0$}&0&n&n& \\
  &0&0&n&n& \\
  &&0&0&n&n \\
  \end{array}
  \]
Consequently, the only relevant vertical differential is
  \begin{align*}
  .df: Ag&\to Ag\ot V\\
  bg&\mapsto
  bw^{-1}yge_x+\sigma(b)xge_y-ba'ge_h.
  \end{align*}
Adding $d(qge_x\wedge e_y)$ one sees that this element is cohomologous to
$(-ba'+ q(\sigma(a')-a') )ge_h$,
where $q\in k[h]$ is such that $(\sigma-w.\Id)(q)=-\sigma(b)$. But then
$b=w\sigma^{-1}(q)-q$, and therefore 
  \[
  .df(bg)=(-w\sigma^{-1}(q)a'+ q\sigma(a') )ge_h=
  (\sigma-w.\Id)(\sigma^{-1}(q)a')ge_h.
  \]
On the other hand, the target of $.df$ has been already shown to be 
isomorphic to $k[h]/(\sigma-w.\Id)(a.k[h])$. Under this isomorphism, the 
cokernel of $.df$ is isomorphic to $k[h]/(\sigma-w.\Id)(a.k[h]+a'.k[h])$. 
Since we have assumed that $(a;a')=1$, the cokernel of $.df$ is zero, 
and by counting dimensions, the kernel of $.df$ also vanishes.
This proves the first part of proposition \ref{prop:coef}; the rest of the
statements thereof follow from similar computations, which we omit.
In particular, theorem \ref{teoG} follows.

\medskip

We would like to observe that in order to prove theorem \ref{teoG}, 
if one assumes that the action of $\Z_g$ is trivial, 
then the
full strength of proposition \ref{prop:coef} is not needed. Indeed, let $g\in
G$. Using proposition \ref{prop:desc} for the cyclic subgroup 
group $C=(g)\subset\Aut_k(A)$
generated by $g$, and the hypothesis on the triviality of the action of the
centralizers in homology, we see that the following relation holds for all
$p\geq0$:
  \begin{align}
  \dim_k\HH_p(A^{C}) 
    & = \dim_k\HH_p(A)^{C} +\sum_{1\leq i<|g|} \dim_kH_p(A,Ag^i)^{C}\label{eq:rel}\\
    & = \dim_k\HH_p(A) +\sum_{1\leq i<|g|} \dim_kH_p(A,Ag^i).\notag 
  \end{align}
Now, in view of lemma \ref{lemma:inv}, the algebra $A^C$ is a GWA, so we already
know its homology. For $p=1$ or $p\geq3$, it vanishes, in particular,
$H_p(A,Ag)=0$; for $p=2$, $\dim_k\HH_2(A^{C})=1=\dim_k\HH_2(A)$, so again we
have $H_2(A,Ag)=0$. Finally, computing $g$-commutators as in the end
of section 5.1, 
it is easy to see that $\dim_kH_0(A,Ag^i)\leq n$ for all $1\leq i<|g|$;
since $\dim_k\HH_0(A^{C})=|g|n-1$ and $\dim_k\HH_0(A)=n-1$, relation
\eqref{eq:rel} forces $\dim_kH_0(A,Ag^i)=n$.

%%%%%%%%%%%%%%%%%%%%%%%%%%%%%%%%%%%%%%%%%%%%%%%%%%%%%%%%%%
\section{Applications}\label{sect:apps}

%%%%%%%%%%%%%%%%%%%%%%%%%%%%%%
\subsection{The usual Weyl algebra}

The results of the previous sections apply to the case when $A=A_1(k)$
($\chr k=0$) and $G$ is an arbitrary finite subgroup of $\Aut_k(A)$, because
in this case the finite order automorphisms of $A$ are always diagonizable,
and $\Z_g$ acts trivially on $H_*(A,Ag)$ for all such $g$. We then recover
the results of \cite{afls}.

%%%%%%%%%%%%%%%%%%%%%%%%%%%%%%
\subsection{Primitive quotients of $\U(\mathfrak{sl}_2)$}

If the polynomial $a$ is of degree two, then $A$ is isomorphic to
one of the maximal primitive quotients of $\U(\mathfrak{sl}_2)$.
In this case, O.~Fleury \cite{odile} has proved that 
the group of automorphisms is isomorphic to the amalgamated product
of $\PSL(2,\CC)$ with a torsion-free group. The
action of $\PSL(2,\CC)$ is the one coming from the adjoint
action of $\SL(2,\CC)$ on $\U(\mathfrak{sl}_2)$. There is then a simple
classification, up to conjugacy, of all finite groups
of automorphisms of $A$: they are the cyclic groups $A_n$, the binary 
dihedral groups
$D_n$, and the binary polyhedral groups $E_6$, $E_7$ and $E_8$;
cf.~\cite{Springer}.

In her thesis, for the regular case, O.~Fleury \cite{odile} has computed, 
case by case the action of the centralizers and in this way
she achieves the computation of $\HH_0(A^G)$. For positive
degrees, following
proposition \ref{prop:desc} one has to compute $H_*(A,Ag)$ and
the action of $\Z_g$ on it. After proposition \ref{prop:coef}
one knows that $H_*(A,Ag)=0$ for $*>0$ and $g\neq 1$, so, for positive
degrees $\HH_*(A^G)=\HH_*(A)^G$. But the only positive and nonzero
degree of $\HH_*(A)$ is $*=2$ and
$\HH_2(A)\cong \HH^0(A)=\Z(A)=k$. Since the duality isomorphism is
$G$-equivariant, the action of $G$ on $\HH_2(A)$ is trivial and
we conclude $\HH_*(A^G)=0$ for $*\neq 0,2$ and $\HH_2(A^G)=k$.
Using the duality one has the cohomology. For the non-regular case,
what we are able to compute is not $\HH_*(A^G)$ but
$\HH_*(A\#G)$. The computation of the action of the centralizers 
is discussed in next sections, because those actions can be described
in general, and also ``explain'' the computations made by O.~Fleury.

\subsection{The Cartan involution\label{section:cartan}}

In the case of $\U({\frak sl}_2)$ there is a special automorphism
(that descends to $B_{\lambda}$) defined by $e\mapsto f$, $f\mapsto e$
and $h\mapsto -h$. For an arbitrary GWA $A$ with defining polinomial
$a(h)$ of degree $n$, there are some particular cases on which a similar
automorphism is defined. In \cite{BavJor} the authors find generators
of the automorphism group of GWAs. It turns out that $\Aut_k(A)$ is
generated
by the torus action 
and exponentials of inner
derivations, and, in the case that there is $\rho\in \CC$
such that
$a(\rho -h)=(-1)^ {n}a(h)$, a generalization
of the Cartan involution --still called $\Om$-- is defined as follows:
\[x\mapsto y \qquad;\qquad y\mapsto (-1)^{n}x\qquad ; \qquad h\mapsto 1+\rho-h\]
Let us call in this section $\G$ the subgroup of $\Aut_k(A)$ generated
by the torus and the exponentials.
When the polynomial is reflective (i.e. $a(\rho-h)=(-1)^na(h)$), 
the group generated by $\G$ and $\Om$ coincides with  $\Aut_k(A)$.
If the polynomial is not reflective, $\G=\Aut_k(A)$.

In \cite{Dixmier} Dixmier shows, in the case of $\U({\frak sl}_2)$, that $\Om$
belongs to $\G$ (and of course this fact descends to the primitive quotients). This
situation corresponds to $\deg(a)=2$. 
In \cite{BavJor} (remark 3.30), the authors ask whether this 
automorphism $\Om$ belongs to the subgroup generated by the torus and exponentials.
We will answer this
question looking at the action of the group of automorphisms
on $\HH_0(A)$, so we begin by recalling the expression of the exponential-type
automorphisms. 

Let $\lambda\in \CC$ and $m\in \NN_0$, the two exponential automorphisms
associated to them are defined as follows:
\begin{align*}
\phi_{m,\lambda}:&=\exp(\lambda \ad( y^m))\\
x&\mapsto x+\sum_{i=1}^n\frac{(-\lambda)^i}{i!}(\ad y^m)^i(x)\\
y&\mapsto y\\
h&\mapsto h+m\lambda y^{m}
\end{align*}
\begin{align*}
\psi_{m,\lambda}:&=\exp(\lambda \ad( x^m))\\
x&\mapsto x\\
y&\mapsto y+\sum_{i=1}^n\frac{\lambda^i}{i!}(\ad x^m)^i(y)\\
h&\mapsto h-m\lambda x^{m}
\end{align*}
We know that $\HH_0(A)$ has $\{1,h,h^2,\dots,h^{n-2}\}$
as a basis. If
one assumes $n>2$ then the action of $\Om$ on $\HH_0(A)$ is not trivial. On the
other hand, since the homogeneous components of weight different
from zero are commutators, the action of $\psi_{m,\lambda}$
and of $\phi_{n,\lambda}$ is trivial on $\HH_0(A)$. 
To see this we consider for example
\[\phi_{m,\lambda}(h^i)= \phi_{m,\lambda}(h)^i= (h+m\lambda y^{m})^i\]
and it is clear that the 0-weight component of $(h+m\lambda y^m)^i$
equals $h^i$.
The action of the torus is also
trivial on $\HH_0(A)$ (it is already trivial on $k[h]$). We conclude
that for $n>2$, the automorphism $\Om$ cannot belong to $\G$.

\subsection{End of the proof of Theorem \ref{teoG}\label{section:end}}

The hypotheses of theorem \ref{teoG} are that in the finite group
$G$ under consideration, every element $g$ is conjugated (in $\Aut_k(A)$)
to an element of the torus.
We will show that the triviality of the action of $\Z_g$ on
$H_*(A,Ag)$ is generically satisfied,
and the fact that the action of $\Z_g$ is trivial or
not depends only on whether $\Om$ belongs to $\Z_g$.

\begin{prop}\label{prop:action}
Let $A$ be a GWA, $G$ a finite subgroup of $\Aut_k(A)$ such that
every $g\in G$ is conjugated to an element of the torus. If
$\Om\notin \Z_g$ then the action of $\Z_g$ on $H_*(A,Ag)$ is trivial.
\end{prop}
\begin{proof}
Let us suppose that $a$ is not reflective and let $g\in G$. 
If $g$ is the identity, then the centralizer of $g$ in $\Aut_k(A)$ is
$\Aut_k(A)$ itself, and
the triviality of the action on $\HH_0$ was explained in the above section.
For $\HH_2(A)\cong \HH^0(A)=\Z(A)=k$ the action is always trivial.
When $*\neq 0,2$, $\HH_*(A)=0$.

Let us now consider $g\neq \Id$. 
Up to conjugation may assume that $g(h)=h$, $g(x)=wx$ and $g(y)=w^{-1}y$, where
$w$ is a root of unity.
After proposition \ref{prop:coef} the
only non-zero homology group is $H_0(A,Ag)$ which has basis
$\{g, hg, h^2g,\dots,h^{n-1}g\}$.

If $g'\in \Aut_k(A)$ commutes with $g$ then it induces an automorphism
$g'|_{A^g}:A^g\to A^g$. This defines a  map $\Z_g\to \Aut_k(A^g)$.
But $A^g=\cl{h, x^{|w|}, y^{|w|} }$ is again a
GWA so one knows the generators of its group of automorphisms.

It is not hard to see that if $a$ is not reflective then the polynomial
$\wt{a}$ associated to $A^g$ is not reflective, neither. %...
In this case $\Aut_k(A)$ is generated by the torus and exponentials
of $\lambda\ad x^{|w|.m}$ 
and $\lambda\ad y^{|w|.m}$. On the other hand, an automorphism $g'$
induces the identity on $A^g$ if and only if it fixes $h$, $x^{|w|}$
and $y^{|w|}$, and by degree considerations it is clear that $g'(x)$ must
be a multiple of $x$ and analogously for $y$. We conclude that the group
of elements (in $\Aut_k(A)$) commuting with $g$ is generated by the torus
and exponentials of $\lambda \ad x^{m.|w|}$ and $\lambda \ad y^{m.|w|}$.

It is clear that the torus acts trivially on the vector space spanned
by $\{g, hg, h^2g,\dots,h^{n-1}g\}$. 
A computation similar to the one done in the end of section \ref{section:cartan}
shows that $\phi_{m,\lambda}(h^ig)=h^ig$ modulo $g$-commutators, 
and analogously for $\psi$.

If the polynomial $a$ is reflective then $\Aut_k(A)$ is generated
by $\G$ and $\Om$. We have just seen the triviality of the action
for the generators of $\G$, and $\Om$ is excluded by hypothesis.
\end{proof}

We now finish the proof of Theorem \ref{teoG}:

We will explain the formula $\dim \HH_0(A^G)=(n-1)+n.a_1+[(n+1)/2].a_2$. 
From the decomposition $\HH_0(A^G)=\bigoplus_{\cl{g}\in
\cl{G}}H_0(A,Ag)^{\Z_g}$, the first ``$n-1$'' comes from the summand
corresponding
to the identity element that contributes with $\dim \HH_0(A)^G=\dim HH_0(A)=n-1$.
The ``$n.a_1$'' comes from the terms corresponding to conjugacy classes
$\cl{g}$ such that $\Z_g$ does not contain $\Om$ because
in this case $\dim H_0(A,Ag)^{\Z_g}=
\dim H_0(A,Ag)=n$.
Finally, the summand ``$[(n+1)/2].a_2$'' corresponds to the conjugacy clases having
$\Om$ in their centralizers, in these cases the dimension of
$(kg\oplus khg\oplus kh^2g\oplus \dots \oplus kh^{n-1}g)^{\Om}$
is the integer part of the half of $n+1$.

\begin{remark}
This result explains, for $n=2$, the case by
case computations of $\HH_0(B_{\lambda}^G)$ made by O.~Fleury 
in \cite{odile}, see also \cite{odile2}.
\end{remark}

%%%%%%%%%%%%%%%%%%%%%%%%%%%%%%

%\subsection{Several variables}
%
%
% In \cite{Bav}, Bavula defines generalized Weyl algebras of arbitrary rank
% $A=D(\sigma,a)$ where $\sigma=(\sigma_1,\dots,\sigma_n)$
% is a set of commuting algebra automorphisms of the ring $D$, and
% $a=(a_1,\dots, a_n)$ is a set of non-zero elements of $\Z(D)$,
% where $\sigma_i(a_j)=a_j$ for all $i\neq j$. The algebra $A$
% is the one obtained from $D$ by adjoining symbols
% $x_1,\dots,x_n,y_1,\dots,y_n$ satisfying the following
% relations:
%  \begin{align*}
%  y_ix_i&=a_i \\
%  x_iy_i&=\sigma_i(a_i) \\
%  [x_i,x_j]&=[y_i,y_j]=[x_i,y_j]=0, \quad\forall i\neq j \\
%  x_i\alpha&=\sigma_i(\alpha)x_i \ \forall \alpha\in D\\
%  y_i\alpha&=\sigma_i^{-1}(\alpha)y_i \ \forall \alpha\in D
%  \end{align*}
% If $D=k[h_1,\dots h_n]$, $a_i\in k[h_i]$ and $\sigma_i(h_i)\in k[h_i]$ for all
% $i=1,\dots,n$, then $A=A(D,\sigma)= A(k[h_1],\sigma_1)\ot\cdots\ot
% A(k[h_n],\sigma_n)$. The computation of $\HH_*(A)$ and $\HH^*(A)$ in this case
% is immediate. For homology, this is due to the K\"unneth formula, and, since
% $A(k[h_i],\sigma_i)$ admits a resolution by projective finitely generated
% bimodules, this is also true for cohomology.
%
% Examples of this kind are the usual Weyl algebra on several variables
% $A_n(k)$, and also
% invariant subalgebras $A_n(k)^g$ where $g$ is any finite order 
% linear automorphism of $A_n(k)$. 

%%%%%%%%%%%%%%%%%%%%%%%%%%%%%%%%%%%%%%%%%%%%%%%%%%%%%%%%%%%%

\end{document}